\documentclass[11pt]{article}
\usepackage[margin=0.5in]{geometry}
\usepackage{times}
\usepackage{hyperref}
\usepackage[numbers]{natbib}
\usepackage{enumitem}
\usepackage{titlesec}
\titleformat{\section}{\large\bfseries}{\thesection.}{0.5em}{}
\titleformat{\subsection}{\normalsize\bfseries}{\thesubsection}{0.5em}{}
\usepackage{graphicx}

\usepackage{lmodern}
\usepackage{titling}
\usepackage{parskip}
\usepackage{xcolor}
\makeatletter
\patchcmd{\@maketitle}
  {\ifx\@empty\@thanks\else\@thanks\fi}
  {\ifx\@empty\@thanks\else\vspace{-1em}{\small\@thanks}\fi}
  {}
  {}
\makeatother

\usepackage{soul}

\title{Modeling Social Systems: \\ Transparency, Reproducibility, and Responsibility}
\author{Maximino Aldana$^{1}$,
Roni Barak Ventura$^{2}$,
Heather Z. Brooks$^{3}$,
Philip S. Chodrow$^{4}$, 
Filipe Georgiou$^{5}$, \\ 
Joseph Johnson$^{6}$, 
Kre\u{s}imir Josić$^{7}$,
Zachary P. Kilpatrick$^{8}$,
Kath Landgren$^{9}$, 
Andrew Nugent$^{10}$, \\
Maurizio Porfiri$^{11}$, 
Nancy Rodriguez$^{8}$, 
Pablo Suárez‑Serrato$^{1,12}$, 
David White$^{13}$,  
Alexander Wiedemann$^{14}$, \\
\& Sam Zhang$^{15}$
}
\begin{document}
\maketitle

\vspace{-1em}
\begin{center}
\footnotesize
$^{1}$Universidad Nacional Autónoma de México,
$^{2}$New Jersey Institute of Technology,
$^{3}$Harvey Mudd College,
$^{4}$Middlebury College,\\
$^{5}$University of Bath,
$^{6}$Carleton College,
$^{7}$University of Houston,
$^{8}$University of Colorado Boulder, 
$^{9}$Stanford University, \\
$^{10}$University of Cambridge,\quad
$^{11}$New York University,
$^{12}$University of California Santa Barbara
$^{13}$Denison University,  \\
$^{14}$Hamline University,
$^{15}$University of Vermont 
\end{center}

\section*{Why Transparent Modeling Matters}

Mathematical models of complex social systems can enrich social scientific theory, inform interventions, and shape policy. 
From voting behavior to economic inequality and urban development, such models influence decisions that affect millions of lives.
Thus, it is especially important to formulate and present them with transparency, reproducibility, and humility.
Modeling in social domains, however, is often uniquely challenging. 
Unlike in physics or engineering, researchers often lack controlled experiments or abundant, clean data. 
Observational data is sparse, noisy, partial, and missing in systematic ways. 
In such an environment, how can we build models that can inform science and decision-making in transparent and responsible ways?

\begin{figure}[t!]
  \centering
  \includegraphics[width=\linewidth]{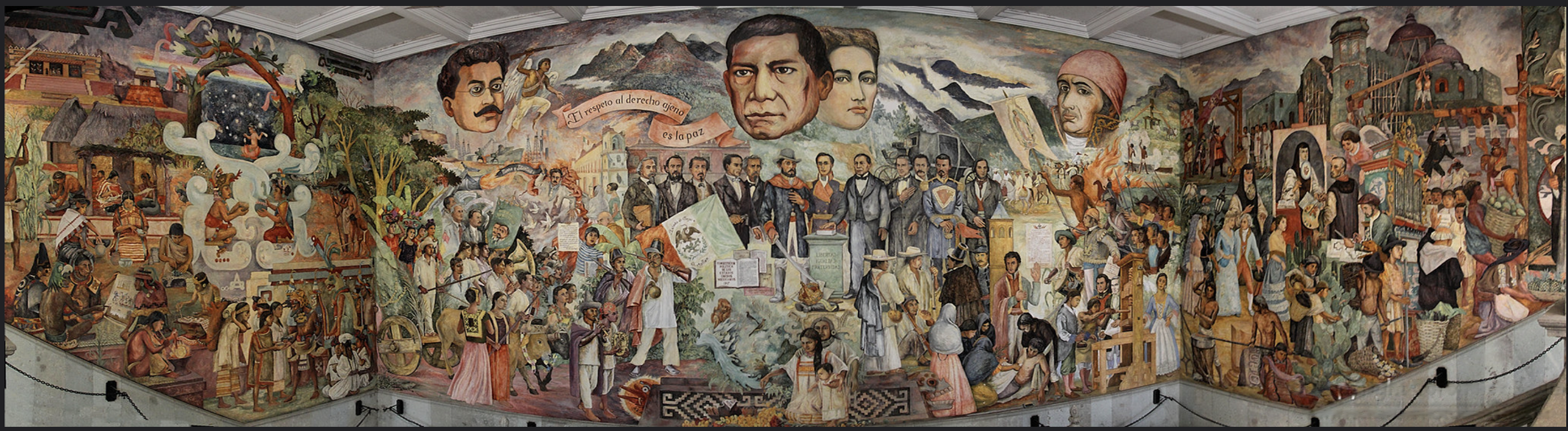}
\caption{
Arturo García Bustos’ mural \textit{Oaxaca en la historia y en el mito}, painted in 1980 in the Palacio de Gobierno in Oaxaca City, a short walk from Casa Matem\'{a}tica Oaxaca. The mural portrays the region’s complex social history from pre-Columbian civilizations through the present. Similarly, mathematical models distill aspects of social systems while omitting others. Unlike art, however, models should make their aims—exploratory, explanatory, or predictive—clear and transparent. (Source: Wikipedia, public domain)
}

  \label{fig:oaxaca_mural}
\end{figure}

This question was the centerpiece of a recent workshop at Casa Matem\'{a}tica Oaxaca titled {\em Collective Social Phenomena: Dynamics and Data.}
The authors of this article---mathematicians, statisticians, engineers, and computational social scientists---spent a week discussing modeling frameworks, ethical challenges, and best practices. 
From these conversations we distilled the following guiding principles for transparent, reproducible, and responsible modeling in social systems (See Fig.~\ref{fig:oaxaca_mural}).
\begin{enumerate}
    \item Be explicit about modeling aims.
    \item Clearly communicate model assumptions and researcher perspectives. 
    \item Match  models  to real-world stakes. 
    \item Quantify and communicate uncertainty.
    \item Share code and data.
    \item Collaborate across disciplines and perspectives.
\end{enumerate}
Below, we expand on each of these principles, drawing examples from our workshop and from broader literature. 
We hope that these principles can serve as conversation-starters and gentle suggestions for quantitative modelers of social systems. 

\section{Why Model?}

``Models are stupid, and we need more of them,'' in the words of Paul Smaldino \citep{smaldino2017}. 
Part of Smaldino's emphasis here is that models may serve purposes other than detailed reconstruction of their target systems, and these purposes may warrant considerable simplifications. 
The appropriateness of a given simplification depends on the purpose the model is intended to serve. 
In the social domain, models can be used to \emph{explore} possible connections between social mechanisms; \emph{explain} observed individual or collective behavior; and \emph{predict} the outcomes of social processes or interventions.
These purposes are distinct from the model's target system or domain of application: one can approach questions of social conformity, voting behavior, economic inequality, or urban development in exploratory, explanatory, or predictive modes.

An example of an \emph{exploratory} model is Schelling's classic model of residential segregation~\citep{schelling1971}, which shows how mild homophilic individual preferences can lead to large-scale segregation patterns. 
Although most scholars attribute racial segregation in the US to long-term discrimination in policy, finance, and urban development, the Schelling model shows that these are sufficient---not necessary---conditions for segregation to emerge. 
The Schelling model can therefore be viewed as an exploration of a counterfactual space of possibility: what \emph{could} happen in a different world.  
Due to its simplicity, the Schelling model is highly interpretable; despite this simplicity, it is also analytically rich. 
The development and analysis of the Schelling model remains an active area of research among mathematicians, quantitative sociologists, and computer scientists.

An \emph{explanatory} model presents (often stylized) descriptions of proposed real-world mechanisms which drive observed social phenomena. 
An example of an explanatory model is the SIR model of infectious disease spread, which describes the spread of disease through a population by dividing individuals into one of three compartments---Susceptible, Infected, or Recovered---and modeling transitions between these compartments. 
Despite its simplicity, the SIR model and its relatives reproduce several stylized features observed in real-world epidemics, such as exponential growth in the early stages of an outbreak and the eventual decline of cases as herd immunity is reached. 
Because the mechanisms in such explanatory models are often highly interpretable, they are often used for scenario analysis, helping policymakers understand the potential impacts of different interventions.

Exploratory and explanatory models can challenge our intuitions and clarify the conditions under which widely accepted narratives break down. 
By making assumptions explicit and specifying causal mechanisms, even simple models can reveal hidden flaws in verbal explanations that align with common sense. 
This capacity to challenge intuition is one of modeling's most important and underappreciated roles. 
As Watts notes, social phenomena often arise from network structure, feedback, or randomness in ways that defy simple narratives~\citep{watts2014}.

In contrast to exploratory and explanatory models, a \emph{predictive} model aims to correctly forecast some facet of a future event based on past data. 
This is often possible without explicitly modeling any underlying mechanisms. 
Unlike simpler mechanistic models, contemporary deep learning models often lack interpretability, making them difficult to scrutinize or trust in policy-relevant settings~\citep{rudin2019}.
Prediction in the social domain is often very challenging. 
Moreover, prediction in the social domain is inherently difficult: the Fragile Families Challenge~\citep{salganik2020} asked 160 teams to predict life outcomes for children based on a common dataset; the resulting models, using widely different approaches, yielded surprisingly similar predictions with surprisingly low predictive accuracy, only narrowly outperforming a simple baseline model. 

Trouble arises when these categories of models are conflated. 
Exploratory models presented as explanations of empirical phenomena can mislead and harm; for example, presenting the Schelling model as an explanatory model of segregation in the US minimizes the role of systemic racism in shaping residential housing patterns. 
Predictive models with strong forecast performance are sometimes taken to obviate simpler explanatory models, echoing claims that abundant data can replace theory altogether~\citep{watts2014}.
On the other hand, relying on interpretable explanatory models for prediction tasks can lead to reduced accuracy in settings where accuracy, rather than interpretability, may be paramount. 
Exploratory, explanatory, and predictive models are valuable in social science when transparently matched to their purpose; usefulness depends less on realism than on clarity of aim.

\section{Communicating Assumptions and Context}

Every model that aims to represent the world contains simplifying assumptions, whether implicit or explicit. 
Transparent modeling requires that these assumptions be clearly communicated.
%We also encourage modelers to reflect on how their training, institutional context, and perspectives shape their work. 
Modelers should be forthcoming about their modeling choices and reflect on how disciplinary and institutional contexts may have shaped them.
For example, the choice of features to include in a predictive model often reflects subjective judgments about which factors are relevant, measurable, and available.

In some cases, researchers may find it helpful to document modeling plans or assumptions in advance. 
While this is not always feasible---particularly in exploratory work---early transparency about motivations, simplifying assumptions, and alternatives considered can clarify the purpose and scope of a model. 
These practices also reduce risks of overfitting or confirmation bias and improve interpretability and accountability~\citep{manski2013}.

Other practices can also help expose implicit choices: documenting modeling decisions, sharing code and intermediate results, or inviting external critique. 
Adversarial collaborations—where researchers with contrasting viewpoints co-design and analyze a model—can clarify how differing modeler assumptions can shape conclusions~\citep{manski2013}.

\section{Match Models to Stakes}

Models of social systems can have consequences for human lives when they are used to inform policy and guide interventions. 
Responsible modeling requires careful consideration of the intended use of the model and the potential consequences of its application.
\begin{itemize}
  \item An exploratory model can cause harm if it is misinterpreted as an empirical explanation. 
  For instance, Schelling's segregation model has sometimes been misused to imply that racial segregation is inevitable or benign, rather than a consequence of racially discriminatory policy choices and structural inequality~\citep{massey1993,vybornova2025}. 
  \item An explanatory model of disease dynamics such as the SIR model can be used to inform public health interventions, which may have positive or negative impacts for population health depending on the model's fidelity to the real world and what consequences of interventions (epidemiological, economic, and social) are considered.
  \item A predictive model like the machine learning models discussed in the Fragile Families Challenge can inform interventions that affect the lives of children, but these interventions may be inappropriate due to model inaccuracy.
  Systemic bias on the basis of demographic attributes like race and gender is also a common attribute of social prediction models. 
\end{itemize}
We therefore encourage modelers to ask: What are the intended uses of this model? 
What are the potential \emph{mis}uses? 
Which groups might use or benefit from this model? 
Which groups might be harmed or excluded by this model? 
Could the way the model is communicated itself cause harm?
Such reflection can guide decisions about formulating, releasing, and communicating models, and help establish norms for responsible modeling.

It is important to remember that one impact of releasing a model is the way in which that model appears in public discourse. 
Divorced from context, mathematical models can mislead; be coopted for political ends; or create a false sense of certainty.
During the COVID-19 pandemic, for example, models played a central role in guiding major decisions. 
Some became lightning rods in the public sphere and were often misrepresented or misunderstood~\citep{ferguson2020}. 
Elsewhere, stylized models of opinion dynamics or economic behavior have appeared in editorials as if they were policy recommendations~\citep{watts2014}.
Other cases are more troubling: For example, machine learning models claiming to predict criminality or sexuality from facial images have been widely discredited, yet their publication has lent credibility to harmful and unfounded ideas~\citep{bergstrom2020}.

While modelers do not always have control over how their work is interpreted or received, they can take simple steps to minimize misinterpretation. 
We can clearly state the aims and simplifying assumptions of our models. 
We can unambiguously label an opinion dynamics model as a toy exploratory model, not a causal account of public opinion formation. 
We can ensure that our published work includes limitations sections, sensitivity analysis, and careful, domain-appropriate language. 
When possible, collaboration with domain experts can help ensure appropriate interpretation and terminology.

Finally, engaging with affected communities can provide perspectives that are otherwise invisible to modelers. Involving those most impacted by a model’s assumptions or applications helps identify potential harms, clarify priorities, and foster accountability. This practice not only strengthens trust but also aligns modeling efforts more closely with the lived realities they aim to represent~\citep{oneil2016}.

\section{Quantify and Communicate Uncertainty}

Whenever possible, modelers should test their models against empirical data and communicate the uncertainty in their inputs and outputs. 
High-stakes modeling in social systems is complicated by challenges common to social data: Low data coverage, systematic missingness, and changing data definitions (such as ethnicity in the US Census) are just a few challenges of empirical modeling. 
Technical tools such as multiple imputation, Bayesian inference, bootstrapping, and confidence intervals can help explicitly represent uncertainty in model inputs and propagate this uncertainty through to model outputs.
Structural modeling or instrumental variable techniques can also support causal inference in the presence of noisy or incomplete data. 
The goal of these approaches is to circumscribe and communicate the uncertainty of models and their alignment with data and verbal theories.
Correct use of these tools typically requires considerable expertise; modelers not trained in these methods may wish to seek out collaborators with expertise in computational social science, econometrics, and causal inference.

\section{Share Code and Data (When Possible)}

It is now common practice to share code and data in software repositories such as GitHub, GitLab, and Bitbucket.
Shared code precisely specifies a model, making it especially helpful for communicating assumptions.
The guidance to publish data can be more nuanced when the data concerns human beings due to questions about privacy, consent, and data ownership.
In some cases in which sharing original data is not possible, modelers can use data synthesizers to construct shareable, synthetic datasets whose population characteristics resemble those of the original data. 

In contexts where models have unambiguous, quantifiable measures of success, as is common in prediction problems, multi-team modeling challenges can be a powerful way to explore many assumptions and approaches in parallel; such approaches have been seen in academic studies like the Fragile Families Challenge~\citep{salganik2020} as well as in public settings like forecasting the COVID-19 pandemic.

\section{Collaborate Across Disciplines and Perspectives}

Modeling may begin as a solo endeavor, but even individual efforts ultimately join a broader conversation. 
Dialogue with collaborators, critics, and communities can improve models, enhance methods, and uncover blind spots, as shown by  initiatives like the Fragile Families Challenge and multi-team replication studies.
Such collaboration is supported by sharing code and data, publishing preprints, and inviting critique from diverse perspectives.

Fortunately, this spirit of openness is increasingly institutionalized. 
The Berkeley Initiative for Transparency in Social Sciences (BITSS), for example, offers tools and training for open science workflows \citep{bitss2020}. 
Reproducibility checklists, data biographies, and adversarial collaborations may all contribute to healthier modeling ecosystems.

\section*{Conclusion: Nurturing a Culture of Responsible Social Modeling}

The \href{https://www.wilmott.com/the-modelers-hippocratic-oath/}{\textit{Hippocratic Oath for Modelers}}, proposed after the 2008 financial crisis, begins: ``I will remember that I didn't make the world, and it doesn't satisfy my equations.'' 
Reasoning about models is not the same as reasoning about the systems they represent, but with appropriate caution—as Smaldino notes—even simple models can illuminate complex behaviors \citep{smaldino2017}.

The toolbox available to modelers grows more potent by the day: generative and predictive artificial intelligence continue to grow at explosive rates; processors for running agent-based simulations become ever-more powerful; and data sets grow more massive, detailed, and networked. 
Against this background, it is essential that modelers reflect on their responsibilities to their stakeholders and to the broader public. 

We have offered six guiding principles for transparent, reproducible, and responsible modeling in social systems: clearly stating modeling aims; communicating model assumptions and disciplinary context; matching models to their real-world stakes; communicating uncertainty; sharing code and data; and collaboration across disciplinary lines. 
We hope that these guiding principles can help researchers at all levels reflect on their own modeling practices and their opportunities to contribute to a culture of responsible modeling of social systems.


\begin{thebibliography}{99}

\bibitem[Smaldino(2017)]{smaldino2017}
Smaldino, P. E. (2017). Models are stupid and we need more of them. \textit{Synthese}, 195(6), 2549–2569.

\bibitem[Ferguson et al.(2020)]{ferguson2020}
Ferguson, N. M., et al. (2020). Impact of non-pharmaceutical interventions (NPIs) to reduce COVID-19 mortality and healthcare demand. \textit{Imperial College COVID-19 Response Team}.

\bibitem[Schelling(1971)]{schelling1971}
Schelling, T. C. (1971). Dynamic models of segregation. \textit{Journal of Mathematical Sociology}, 1(2), 143–186.

\bibitem[Rudin(2019)]{rudin2019}
Rudin, C. (2019). Stop explaining black box machine learning models for high stakes decisions and use interpretable models instead. \textit{Nature Machine Intelligence}, 1(5), 206–215.


\bibitem[Salganik et al.(2020)]{salganik2020}
Salganik, M. J., et al. (2020). Measuring the predictability of life outcomes with a scientific mass collaboration. \textit{PNAS}, 117(15), 8398–8403.

\bibitem[Jumper et al.(2021)]{jumper2021}
Jumper, J., et al. (2021). Highly accurate protein structure prediction with AlphaFold. \textit{Nature}, 596, 583–589.



\bibitem[Manski(2013)]{manski2013}
Manski, C. F. (2013). \textit{Public Policy in an Uncertain World: Analysis and Decisions}. Harvard University Press.

\bibitem[Massey \& Denton(1993)]{massey1993}
Massey, D. S., \& Denton, N. A. (1993). \textit{American Apartheid: Segregation and the Making of the Underclass}. Harvard University Press.

\bibitem[Vybornova \& Verma(2025)]{vybornova2025}
Vybornova, A., \& Verma, T. (2025). You Don’t Have to Live Next to Me: Towards Demobilizing Individualistic Bias in Computational Approaches to Urban Segregation. \textit{arXiv preprint arXiv:2505.01830}.



\bibitem[BITSS(2020)]{bitss2020}
Berkeley Initiative for Transparency in the Social Sciences (BITSS). (2020). Retrieved from \url{https://www.bitss.org}

\bibitem[Watts(2014)]{watts2014}
Watts, D. J. (2014). Common sense and sociological explanations. \textit{American Journal of Sociology}, 120(2), 313–351.

\bibitem[Bergstrom \& West(2020)]{bergstrom2020}
Bergstrom, C. T., \& West, J. D. (2020). \textit{Calling Bullshit: The Art of Skepticism in a Data-Driven World}. Random House.

\bibitem[O'Neil(2016)]{oneil2016}
O'Neil, C. (2016). \textit{Weapons of Math Destruction: How Big Data Increases Inequality and Threatens Democracy}. Crown Publishing Group.


\end{thebibliography}
\end{document}